\documentclass[a4paper]{amsart}
\usepackage[T1]{fontenc}
\usepackage{amssymb,enumitem,mathrsfs,mathtools}
\usepackage[all]{xy}
\usepackage[ngerman,british]{babel}
\usepackage{lmodern}
\usepackage{microtype}

\title{Coarse and equivariant co-assembly maps}

\author{Heath Emerson}
\email{hemerson@math.uvic.ca}
\address{Department of Mathematics and Statistics\\
         University of Victoria\\
         PO BOX 3045 STN CSC \\
         Victoria, B.C.\\
         Canada \\
         V8W 3P4}

\author{Ralf Meyer}
\email{rameyer@uni-math.gwdg.de}
\address{\selectlanguage{ngerman}Mathematisches Institut\\
         Georg-August-Universit\"at G\"ottingen\\
         Bunsenstra\ss e 3--5\\
         37073 G\"ottingen\\
         Deutschland}

\usepackage[pdftitle={Coarse and equivariant co-assembly maps},
  pdfauthor={Heath Emerson and Ralf Meyer},
  pdfsubject={Mathematics; MSC 19K35, 46L80}
]{hyperref}
\usepackage[lite]{amsrefs}
\newcommand*{\MRref}[2]{\linebreak[0] \href{http://www.ams.org/mathscinet-getitem?mr=#1}{MR \textbf{#1}}}

\theoremstyle{plain}
\newtheorem{theorem}{Theorem}
\newtheorem{lemma}[theorem]{Lemma}
\newtheorem{proposition}[theorem]{Proposition}
\newtheorem{corollary}[theorem]{Corollary}
\newtheorem{definition}[theorem]{Definition}
\theoremstyle{remark}

\newtheorem{example}[theorem]{Example}

\newcommand*{\C}{\mathbb C}
\newcommand*{\Z}{\mathbb Z}
\newcommand*{\Ztwo}{\mathbb Z/2}
\newcommand*{\R}{\mathbb R}
\newcommand*{\N}{\mathbb N}

\newcommand*{\vvr}{\mathfrak{\bar c}^\mathfrak{red}}
\newcommand*{\bvr}{\mathfrak c^\mathfrak{red}}

\newcommand*{\bbbarred}{\bar{\mathfrak B}^\mathfrak{red}}
\newcommand*{\bbred}{\mathfrak B^\mathfrak{red}}
\newcommand*{\bbbar}{\bar{\mathfrak B}}

\newcommand*{\ADir}{\mathsf P}
\newcommand*{\Dirac}{\mathsf D}

\newcommand*{\KK}{\textup{KK}}
\newcommand*{\RKK}{\textup{RKK}}
\newcommand*{\K}{\textup K}
\newcommand*{\Ktop}{\textup K^\textup{top}}
\newcommand*{\KX}{\textup{KX}}

\newcommand*{\point}{\textup{point}}
\newcommand*{\Jet}{\textup J}
\newcommand*{\Spin}{\textup{Spin}}
\newcommand*{\ID}{\textup{id}}

\DeclareMathOperator{\Ind}{Ind}
\DeclareMathOperator{\Rep}{Rep}
\DeclareMathOperator{\supp}{supp}
\DeclareMathOperator{\Diff}{Diff}

\newcommand*{\Hilm}{\mathcal E}
\newcommand*{\EG}{\mathcal EG}

\newcommand*{\Rips}{\mathscr P}
\newcommand*{\Subrips}{\textup P}
\newcommand*{\X}{\mathscr X}

\newcommand*{\Mult}{\mathcal M}
\newcommand*{\Cmax}{C^*_\textup{max}}

\newcommand*{\Comp}{\mathbb K}

\newcommand*{\nb}{\nobreakdash-}

\newcommand*{\abs}[1]{\lvert#1\rvert}
\newcommand*{\norm}[1]{\lVert#1\rVert}

\newcommand*{\maxotimes}{\mathbin{\otimes_\textup{max}}}
\newcommand*{\rcross}{\mathbin{\rtimes_\textup r}}
\newcommand*{\cross}{\mathbin\rtimes}
\newcommand*{\defeq}{\mathrel{\vcentcolon=}}
\newcommand*{\into}{\rightarrowtail}
\newcommand*{\prto}{\twoheadrightarrow}
\newcommand*{\blank}{\text\textvisiblespace}

\begin{document}

\begin{abstract}
  We study an equivariant co-assembly map that is dual to the usual
  Baum--Connes assembly map and closely related to coarse geometry, equivariant
  Kasparov theory, and the existence of dual Dirac morphisms.  As
  applications, we prove the existence of dual Dirac morphisms for groups with
  suitable compactifications, that is, satisfying the Carlsson--Pedersen
  condition, and we study a \(\K\)\nb-theoretic counterpart to the proper
  Lipschitz cohomology of Connes, Gromov and Moscovici.
\end{abstract}

%% \begin{classification}
%%   Primary 19K35; Secondary 46L80.
%% \end{classification}

%% \begin{keywords}
%%   K-theory, coarse geometry, Baum--Connes conjecture, assembly map,
%%   compactification.
%% \end{keywords}



\thanks{This research was partially carried out at the
  {\selectlanguage{ngerman}Westf\"alische Wilhelms-Universit\"at M\"unster}
  and supported by the EU-Network \emph{Quantum Spaces and Noncommutative
    Geometry} (Contract HPRN-CT-2002-00280) and the
  \emph{\selectlanguage{ngerman}Deutsche Forschungsgemeinschaft} (SFB 478).}

\maketitle

\section{Introduction}
\label{sec:intro}

This is a sequel to the articles \cites{Emerson-Meyer:Dualizing,Emerson-Meyer:Descent}, which
deal with a coarse co-assembly map that is dual to the usual coarse assembly
map.  Here we study an equivariant co-assembly map that is dual to the
Baum--Connes assembly map for a group~\(G\).

A rather obvious choice for such a dual map is the map
\begin{equation}
  \label{eq:dumb_co-assembly}%
  p_\EG^*\colon \KK^G_*(\C,\C)\to\RKK^G_*(\EG;\C,\C)
\end{equation}
induced by the projection \(p_\EG\colon \EG\to\point\).  This map and its
application to the Novikov conjecture go back to Kasparov (\cite{Kasparov:Novikov}).
Nevertheless, \eqref{eq:dumb_co-assembly} is not quite the map that we
consider here.  Our map is closely related to the coarse co-assembly map
of~\cite{Emerson-Meyer:Dualizing}.  It is an isomorphism if and only if the
Dirac-dual-Dirac method applies to~\(G\).  Hence there are many cases\mdash
groups with \(\gamma\neq1\)\mdash where our co-assembly map is an
isomorphism and~\eqref{eq:dumb_co-assembly} is not.

Most of our results only work if the group~\(G\) is (almost) totally
disconnected and has a \(G\)\nb-compact universal proper \(G\)\nb-space
\(\EG\).  We impose this assumption throughout the introduction.

First we briefly recall some of the main ideas
of~\cites{Emerson-Meyer:Dualizing,Emerson-Meyer:Descent}.  The new ingredient in the coarse
co-assembly map is the \emph{reduced stable Higson corona} \(\bvr(X)\) of a
coarse space~\(X\).  Its definition resembles that of the usual Higson
corona, but its \(\K\)\nb-theory behaves much better.  The coarse
co-assembly map is a map
\begin{equation}
  \label{eq:coarse_co-assembly}
  \mu\colon \K_{*+1}\bigl(\bvr(X)\bigr) \to \KX^*(X),
\end{equation}
where \(\KX^*(X)\) is a coarse invariant of~\(X\) that agrees with \(\K^*(X)\)
if~\(X\) is uniformly contractible.

If~\(\abs{G}\) is the coarse space underlying a group~\(G\), then there is a commuting
diagram
\begin{equation}
  \label{eq:non-equivariant_co-assembly}
  \begin{gathered}
  \xymatrix@C+2em{
    \KK^G_*\bigl(\C,C_0(G)\bigr) \ar[r]^-{p_\EG^*} \ar@{<->}[d]^{\cong} &
    \RKK^G_*\bigl(\EG;\C,C_0(G)\bigr) \ar@{<->}[d]^{\cong} \\
    \K_{*+1}\bigl(\bvr(\abs{G})\bigr) \ar[r]^-{\mu} &
    \KX^*(\abs{G}).
  }
  \end{gathered}
\end{equation}
In this situation, \(\KX^*(\abs{G})\cong \K^*(\EG)\) 
because~\(\abs{G}\) is coarsely
equivalent to~\(\EG\), which is uniformly contractible.  The commuting
diagram~\eqref{eq:non-equivariant_co-assembly}, coupled with the 
reformulation of the Baum Connes assembly map in 
\cite{Meyer-Nest:BC}, is the source of the
relationship between the coarse co-assembly map and the Dirac-dual-Dirac
method mentioned above. 

If~\(G\) is a torsion-free discrete group with finite classifying space
\(BG\), the coarse co-assembly map is an isomorphism if and only if the
Dirac-dual-Dirac method applies to~\(G\).  A similar result for
groups with torsion is available, but this requires
 working equivariantly with respect to compact
subgroups of~\(G\).

In this article, we work equivariantly with respect to the whole group~\(G\).  
The action
of~\(G\) on its underlying coarse space~\(\abs{G}\) by isometries induces an
action on \(\bvr(G)\).  We consider a \(G\)\nb-equivariant analogue
\begin{equation}
  \label{eq:equivariant_co-assembly}%
  \mu\colon \Ktop_{*+1}\bigl(G,\bvr(\abs{G})\bigr) \to \K_*(C_0(\EG)\cross G)
\end{equation}
of the coarse co-assembly map~\eqref{eq:coarse_co-assembly}; here
\(\Ktop_*(G,A)\) denotes the domain of the Baum--Connes assembly map for~\(G\)
with coefficients~\(A\).  We avoid \(\K_*(\bvr(X)\cross G)\) and
\(\K_*(\bvr(X)\rcross G)\) because we can say nothing about these two groups.
In contrast, the group \(\Ktop_*\bigl(G,\bvr(\abs{G}\bigr)\) is much more
manageable.  The only \emph{analytical} difficulties in this group come from
coarse geometry.

There is a commuting diagram similar to~\eqref{eq:non-equivariant_co-assembly}
that relates~\eqref{eq:equivariant_co-assembly} to equivariant Kasparov
theory.  To formulate this, we need some results of~\cite{Meyer-Nest:BC}.  There
is a certain \(G\)\nb-\(C^*\)-algebra~\(\ADir\) and a class
\(\Dirac\in\KK^G(\ADir,\C)\) called \emph{Dirac morphism} such that the
Baum--Connes assembly map for~\(G\) is equivalent to the map
\[
\K_*\bigl( (A\otimes \ADir) \rcross G) \to \K_*(A\rcross G)
\]
induced by Kasparov product with~\(\Dirac\).  The Baum--Connes conjecture holds
for~\(G\) with coefficients in \(\ADir\otimes A\) for any~\(A\).  The Dirac
morphism is a \emph{weak equivalence}, that is, its image in
\(\KK^H(\ADir,\C)\) is invertible for each compact subgroup~\(H\) of~\(G\).

The existence of the Dirac morphism allows us to localise the (triangulated)
category \(\KK^G\) at the multiplicative system of weak equivalences.  The
functor from \(\KK^G\) to its localisation turns out to be equivalent to the
map
\[
p_\EG^*\colon \KK^G(A,B)\to\RKK^G(\EG; A,B).
\]
One of the main results of this paper is a commuting diagram
\begin{equation}
  \label{eq:maintheorem}%
  \begin{gathered}
    \xymatrix{
      \Ktop_{*+1}\bigl(G,\bvr(\abs{G})\bigr) \ar[r] \ar@{<->}[d]^{\cong} &
      \K^*(\EG) \ar@{<->}[d]^{\cong} \\
      \KK^G_*(\C,\ADir) \ar[r]^-{p_\EG^*} &
      \RKK^G_*(\EG;\C,\ADir).
    }
  \end{gathered}
\end{equation}
In other words, the equivariant coarse
co-assembly~\eqref{eq:equivariant_co-assembly} is equivalent to the map
\[
p_\EG^*\colon \KK^G_*(\C,\ADir) \to \RKK^G_*(\EG;\C,\ADir).
\]
This map is our proposal for a dual to the Baum--Connes assembly map.

We should justify why we prefer the map~\eqref{eq:equivariant_co-assembly}
over~\eqref{eq:dumb_co-assembly}.  Both maps have isomorphic targets:
\[
\RKK^G_*(\EG;\C,\ADir) \cong \RKK^G_*(\EG;\C,\C) \cong \K_*(C_0(\EG)\cross G).
\]
Even in the usual Baum--Connes assembly map, the analytical side involves a
choice between full and reduced group \(C^*\)\nb-algebras and crossed
products.  Even though the full group \(C^*\)\nb-algebra has better
functoriality properties and is sometimes preferred because it gives
potentially finer invariants, the reduced one is used because its
\(\K\)\nb-theory is closer to \(\Ktop_*(G)\). In formulating a dual version
of the assembly map, we are faced with a similar situation.  Namely, the
topological object that is dual to \(\Ktop_*(G)\) is \(\RKK^G_*(\EG;\C,\C)\).
For the analytical side, we have some choices; we prefer \(\KK^G_*(\C,\ADir)\)
over \(\KK^G_*(\C,\C)\) because the resulting co-assembly map is an
isomorphism in more cases.

Of course, we must check that this choice is analytical enough to be useful
for applications. The most important of these is the Novikov conjecture.
Elements of \(\RKK^G_*(\EG;\C,\C)\) yield maps \(\Ktop_*(G)\to\Z\), which are
analogous to higher signatures.  In particular,
\eqref{eq:equivariant_co-assembly} gives rise to such objects.  The maps
\(\Ktop_*(G)\to\Z\) that come from a class in the range
of~\eqref{eq:dumb_co-assembly} are known to yield homotopy invariants for
manifolds because there is a pairing between \(\KK^G_*(\C,\C)\) and
\(\K_*(\Cmax G)\) (see~\cite{Ferry-Ranicki-Rosenberg:Novikov}).  But since
\eqref{eq:equivariant_co-assembly} factors
through~\eqref{eq:dumb_co-assembly}, the former also produces homotopy
invariants.

In particular, surjectivity of \eqref{eq:equivariant_co-assembly} implies the
Novikov conjecture for~\(G\). More is true: since \(\KK^G(\C,\ADir)\) is the
home of a \emph{dual-Dirac morphism}, \eqref{eq:maintheorem} yields that~\(G\)
has a dual-Dirac morphism and hence a \emph{\(\gamma\)\nb-element} if and
only if~\eqref{eq:equivariant_co-assembly} is an isomorphism.  This
observation can be used to give an alternative proof of the main result
of~\cite{Emerson-Meyer:Descent}.

We call elements in the range of~\eqref{eq:equivariant_co-assembly}
\emph{boundary classes}.  These automatically form a graded ideal in the
\(\Ztwo\)-graded unital ring \(\RKK^G_*(\EG;\C,\C)\).  In contrast, the
range of the unital ring homomorphism~\eqref{eq:dumb_co-assembly} need not be
an ideal because it always contains the unit element of
\(\RKK^G_*(\EG;\C,\C)\).  We describe two important constructions of boundary
classes, which are related to compactifications and to the proper Lipschitz
cohomology of~\(G\) studied in \cites{Connes-Gromov-Moscovici:Lipschitz, Dranishnikov:Lipschitz}.

Let \(\EG\subseteq Z\) be a \(G\)\nb-equivariant compactification of \(\EG\)
that is compatible with the coarse structure in a suitable sense.  Since there
is a map
\[
\Ktop_*\bigl(G,C(Z\setminus \EG)\bigr) \to \Ktop_*\bigl(G,\bvr(\EG)\bigr),
\]
we get boundary classes from the boundary \(Z\setminus \EG\).  This
construction also shows that~\(G\) has a dual-Dirac morphism if it satisfies
the Carlsson--Pedersen condition.  This improves upon a result of Nigel Higson
(\cite{Higson:Bivariant}), which shows split injectivity of the
Baum--Connes assembly map with coefficients under the same assumptions.

Although we have discussed only \(\KK^G_*(\C,\ADir)\) so far, our main
technical result is more general and can also be used to construct elements in
Kasparov groups of the form \(\KK^G_*\bigl(\C,C_0(X)\bigr)\) for suitable
\(G\)\nb-spaces~\(X\).  If~\(X\) is a proper \(G\)\nb-space, then we can
use such classes to construct boundary classes in \(\RKK^G_*(\EG;\C,\C)\).
This provides a \(\K\)\nb-theoretic counterpart of the proper Lipschitz
cohomology of~\(G\) defined by Connes, Gromov, and Moscovici in~\cite{Connes-Gromov-Moscovici:Lipschitz}.
Our approach clarifies the geometric parts of several constructions
in~\cite{Connes-Gromov-Moscovici:Lipschitz}; thus we substantially simplify the proof of the homotopy
invariance of Gelfand--Fuchs cohomology classes in~\cite{Connes-Gromov-Moscovici:Lipschitz}.

\section{Preliminaries}
\label{sec:preliminaries}

\subsection{Dirac-dual-Dirac method and Baum--Connes conjecture}
\label{sec:Dirac_and_BC}

First, we recall the \emph{Dirac-dual-Dirac method} of Kasparov and its
reformulation in~\cite{Meyer-Nest:BC}.  This is a technique for proving
injectivity of the \emph{Baum--Connes assembly map}
\begin{equation}
  \label{eq:BC_assembly}
  \mu\colon \Ktop_*(G,B) \to \K_*\bigl(C^*_r(G,B)\bigr),
\end{equation}
where~\(G\) is a locally compact group and~\(B\) is a \(C^*\)\nb-algebra
with a strongly continuous action of~\(G\) or, briefly, a
\emph{\(G\)\nb-\(C^*\)-algebra}.

This method requires a proper \(G\)-\(C^*\)-algebra~\(A\) and classes
\[
d\in\KK^G(A,\C),
\qquad
\eta\in\KK^G(\C,A),
\qquad
\gamma \defeq \eta \otimes_A d \in \KK^G(\C,\C),
\]
such that \(p_\EG^*(\gamma) = 1_\C\) in \(\RKK^G(\EG;\C,\C)\).  If these data
exist, then the Baum--Connes assembly map~\eqref{eq:BC_assembly} is injective
for all~\(B\).  If, in addition, \(\gamma=1_\C\) in \(\KK^G(\C,\C)\), then the
Baum--Connes assembly map is invertible for all~\(B\), so that~\(G\) satisfies
the Baum--Connes conjecture with arbitrary coefficients.

Let \(A\) and~\(B\) be \(G\)\nb-\(C^*\)-algebras.  An element \(f\in
\KK^G(A,B)\) is called a \emph{weak equivalence} in~\cite{Meyer-Nest:BC} if its
image in \(\KK^H(A,B)\) is invertible for each compact subgroup~\(H\)
of~\(G\).

The following theorem contains some of the main results of~\cite{Meyer-Nest:BC}.

\begin{theorem}
  \label{meyernest}
  Let~\(G\) be a locally compact group.  Then there is a
  \(G\)\nb-\(C^*\)-algebra~\(\ADir\) and a class
  \(\Dirac\in\KK^G(\ADir,\C)\) called \emph{Dirac morphism} such that
  \begin{enumerate}[label=\textup{(\alph*)}]
  \item \(\Dirac\) is a weak equivalence;

  \item the Baum--Connes conjecture holds with coefficients in
    \(A\otimes\ADir\) for any~\(A\);

  \item the assembly map~\eqref{eq:BC_assembly} is equivalent to the map
    \[
    \Dirac_*\colon \K_*(A\otimes\ADir\rcross G) \to \K_*(A\rcross G);
    \]

  \item the Dirac-dual-Dirac method applies to~\(G\) if and only if there is a
    class \(\eta\in\KK^G(\C,\ADir)\) with \(\eta\otimes_\C D = 1_A\), if and
    only if the map
    \begin{equation}
      \label{diracdualdiracmethod}
      \Dirac^*\colon \KK^G_*(\C,\ADir) \to \KK^G_*(\ADir,\ADir),
      \qquad
      x\mapsto D\otimes x,
    \end{equation}
    is an isomorphism.
  \end{enumerate}
\end{theorem}

Whereas~\cite{Emerson-Meyer:Descent} studies the invertibility
of~\eqref{diracdualdiracmethod} by relating it
to~\eqref{eq:coarse_co-assembly}, here we are going to study the
map~\eqref{diracdualdiracmethod} itself.

It is shown in~\cite{Meyer-Nest:BC} that the localisation of the category
\(\KK^G\) at the weak equivalences is isomorphic to the category
\(\RKK^G(\EG)\) whose morphism spaces are the groups \(\RKK^G(\EG;A,B)\) as
defined by Kasparov in~\cite{Kasparov:Novikov}.  This statement is equivalent to the
existence of a Poincar\'e duality isomorphism
\begin{equation}
  \label{eq:PD_EG}
  \KK^G_*(A\otimes\ADir,B) \cong \RKK^G_*(\EG;A,B)
\end{equation}
for all \(G\)-\(C^*\)-algebras \(A\) and~\(B\) (this notion of duality is
analysed in~\cite{Emerson-Meyer:Euler}).  The canonical functor from \(\KK^G\)
to the localisation becomes the obvious functor
\[
p_\EG^*\colon \KK^G(A,B) \to \RKK^G(\EG; A,B).
\]

Since~\(\Dirac\) is a weak equivalence, \(p_\EG^*(\Dirac)\) is invertible.
Hence the maps in the following commuting square are isomorphisms for all
\(G\)\nb-\(C^*\)-algebras \(A\) and~\(B\):
\[\xymatrix@C=4em{
  \RKK^G_*(\EG; A, B\otimes \ADir) \ar[r]_{\cong}^{\Dirac_*} \ar[d]_{\cong}^{\Dirac^*} &
  \RKK^G_*(\EG; A,B) \ar[d]_{\cong}^{\Dirac^*} \\
  \RKK^G_*(\EG; A\otimes\ADir,B\otimes\ADir) \ar[r]_{\cong}^{\Dirac_*} &
  \RKK^G_*(\EG; A\otimes\ADir,B).
}
\]
Together with~\eqref{eq:PD_EG} this implies
\[
\KK^G_*(A\otimes\ADir,B)\cong\KK^G_*(A\otimes\ADir, B\otimes\ADir).
\]
In the following, it will be useful to turn the isomorphism
\[
\Ktop_*(G,A) \cong \K_*\bigl((A\otimes\ADir)\rcross G\bigr)
\]
in Theorem \ref{meyernest}.(c) into a definition.

\subsection{Group actions on coarse spaces}
\label{sec:acts_coarse}

Let~\(G\) be a locally compact group and let~\(X\) be a right \(G\)\nb-space
and a coarse space.  We always assume that~\(G\) acts continuously and
coarsely on~\(X\), that is, the set \(\{(xg, yg)\mid g \in K, (x,y) \in E\}\)
is an entourage for any compact subset~\(K\) of~\(G\) and any entourage~\(E\)
of~\(X\).

\begin{definition}
  \label{typesofactions}
  We say that~\(G\) \emph{acts by translations} on~\(X\) if \(\{(x,gx)\mid
  x\in X,\ g\in K\}\) is an entourage for all compact subsets \(K\subseteq
  G\).  We say that \(G\) acts by \emph{isometries} if every entourage
  of~\(X\) is contained in a \(G\)\nb-invariant entourage.
\end{definition}

\begin{example}
  \label{groupsascoarsespaces}
  Let~\(G\) be a locally compact group.  Then~\(G\) has a unique coarse
  structure for which the right translation action is isometric; the
  corresponding coarse space is denoted~\(\abs{G}\).  The generating
  entourages are of the form
  \[
  \bigcup_{g\in G} Kg\times Kg =
  \{(xg,yg) \mid g\in G, x,y\in K\}
  \]
  for compact subsets~\(K\) of~\(G\).  The \emph{left} translation action is
  an action by translations for this coarse structure.
\end{example}

\begin{example}
  \label{exa:coarse_on_proper_G-space}%
  More generally, any proper, \(G\)\nb-compact \(G\)\nb-space~\(X\)
  carries a unique coarse structure for which~\(G\) acts isometrically; its
  entourages are defined as in Example~\ref{groupsascoarsespaces}.  With this
  coarse structure, the orbit map \(G\to X\), \(g\mapsto g\cdot x\), is a
  coarse equivalence for any choice of \(x\in X\).  If the
  \(G\)\nb-compactness assumption is omitted, the result is a
  \(\sigma\)\nb-coarse space.  We always equip a proper \(G\)\nb-space
  with this additional structure.
\end{example}

\subsection{The stable Higson corona}
\label{sec:Higson_corona}

We next recall the definition of the \emph{stable Higson corona} of a coarse
space \(X\) from \cites{Emerson-Meyer:Dualizing, Emerson-Meyer:Descent}.  Let~\(D\) be a
\(C^*\)\nb-algebra.

Let \(\Mult(D\otimes\Comp)\) be the multiplier algebra of \(D\otimes\Comp\),
and let \(\bbbarred(X,D)\) be the \(C^*\)\nb-algebra of norm-continuous,
bounded functions \(f\colon X\to\Mult(D\otimes\Comp)\) for which
\(f(x)-f(y)\in D\otimes\Comp\) for all \(x,y\in X\).  We also let
\[
\bbred(X,D) \defeq \bbbarred(X,D)/C_0(X,D\otimes\Comp).
\]

\begin{definition}
  A function \(f\in\bbbarred(X,D)\) has \emph{vanishing variation} if the
  function \(E \ni (x,y)\mapsto \norm{f(x)-f(y)}\) vanishes at~\(\infty\) for
  any closed entourage \(E\subseteq X\times X\).
\end{definition}

The \emph{reduced stable Higson compactification} of~\(X\) with
coefficients~\(D\) is the subalgebra \(\vvr(X,D) \subseteq \bbbarred(X,D)\) of
vanishing variation functions.  The quotient
\[
\bvr(X,D) \defeq \vvr(X,D)/C_0(X,D\otimes\Comp) \subseteq \bbred(X,D)
\]
is called \emph{reduced stable Higson corona} of~\(X\).  This defines a
functor on the coarse category of coarse spaces: a coarse map \(f\colon X\to
X'\) induces a map \(\bvr(X',D)\to\bvr(X,D)\), and two maps \(X\to X'\) induce
the same map \(\bvr(X',D)\to\bvr(X,D)\) if they are close.  Hence a coarse
equivalence \(X\to X'\) induces an isomorphism \(\bvr(X',D)\cong\bvr(X,D)\).

For some technical purposes, we must allow unions \(\X=\bigcup X_n\) of coarse
spaces such that the embeddings \(X_n\to X_{n+1}\) are coarse equivalences;
such spaces are called \emph{\(\sigma\)\nb-coarse spaces}.  The main example
is the \emph{Rips complex} \(\Rips(X)\) of a coarse space~\(X\), which is used
to define its coarse \(\K\)\nb-theory.  More generally, if~\(X\) is a proper
but not \(G\)\nb-compact \(G\)\nb-space, then~\(X\) may be endowed with
the structure of a \(\sigma\)\nb-coarse space.  For coarse spaces of the
form~\(\abs{G}\) for a locally compact group~\(G\) with a \(G\)\nb-compact
universal proper \(G\)\nb-space \(\EG\), we may use \(\EG\) instead of
\(\Rips(X)\) because \(\EG\) is coarsely equivalent to~\(G\) and uniformly
contractible.  Therefore, we do not need \(\sigma\)\nb-coarse spaces much;
they only occur in Lemma~\ref{lem:coarse_gives_classes}.

It is straightforward to extend the definitions of \(\vvr(X,D)\) and
\(\bvr(X,D)\) to \(\sigma\)\nb-coarse spaces (see \cites{Emerson-Meyer:Dualizing,
  Emerson-Meyer:Descent}).  Since we do not use this generalisation much, we omit
details on this.

Let~\(H\) be a locally compact group that acts coarsely and properly on~\(X\).
It is crucial for us to allow non-compact groups here, whereas
\cite{Emerson-Meyer:Descent} mainly needs equivariance for compact groups.  Let~\(D\)
be an \(H\)\nb-\(C^*\)-algebra, and let \(\Comp_H\defeq \Comp
(\ell^2\N\otimes L^2H)\).  Then~\(H\) acts on \(\bbbarred(X,D\otimes
\Comp_H)\) by
\[
(h\cdot f)(x)\defeq h\cdot \bigl(f(xh)\bigr),
\]
where we use the obvious action of~\(H\) on \(D\otimes\Comp_H\) and its
multiplier algebra.  The action of~\(H\) on \(\bbbarred(X, D\otimes\Comp_H)\)
need not be continuous; we let \(\bbbarred_H(X,D)\) be the subalgebra of
\(H\)\nb-continuous elements in \(\bbbarred(X, D\otimes\Comp_H)\).  We let
\(\vvr_H(X,D)\) be the subalgebra of vanishing variation functions in
\(\bbbarred_H(X,D)\).  Both algebras contain \(C_0(X, D\otimes\Comp_H)\) as an
ideal.  The corresponding quotients are denoted by \(\bbred_H(X,D)\) and
\(\bvr_H(X,D)\).  By construction, we have a natural morphism of extensions of
\(H\)\nb-\(C^*\)-algebras
\begin{equation}
  \label{eq:forget-control}%
  \begin{gathered}
    \xymatrix{
      C_0(X, D\otimes\Comp_H)\ \ar@{>->}[r] \ar@{=}[d] &
      \vvr_H(X,D) \ar@{->>}[r] \ar[d]^{\subseteq} &
      \bvr_H(X,D) \ar[d]^{\subseteq} \\
      C_0(X, D\otimes\Comp_H)\ \ar@{>->}[r] &
      \bbbarred_H(X,D) \ar@{->>}[r] &
      \bbred_H(X,D).
    }
  \end{gathered}
\end{equation}

Concerning the extension of this construction to \(\sigma\)\nb-coarse
spaces, we only mention one technical subtlety.  We must extend the functor
\(\Ktop_*(H,\blank)\) from \(C^*\)\nb-algebras to
\(\sigma\)\nb-\(H\)-\(C^*\)\nb-algebras.  Here we use the definition
\begin{equation}
  \label{vanilla}
  \Ktop_*(H,A) \cong \K_*\bigl((A\otimes\ADir)\rcross H\bigr),
\end{equation}
where \(\Dirac\in\KK^H(\ADir,\C)\) is a Dirac morphism
for~\(H\).  The more traditional definition as a colimit of
\(\KK^G_*(C_0(X),A)\), where \(X\subseteq\EG\) is
\(G\)\nb-compact, yields a wrong result if~\(A\) is a
\(\sigma\)\nb-\(H\)-\(C^*\)-algebra because colimits and
limits do not commute.

Let~\(H\) be a locally compact group, let~\(X\) be a coarse space with an
isometric, continuous, proper action of~\(H\), and let~\(D\) be an
\(H\)\nb-\(C^*\)-algebra.  The \(H\)\nb-equivariant coarse
\(\K\)\nb-theory \(\KX^*_H(X,D)\) of~\(X\) with coefficients in~\(D\) is
defined in~\cite{Emerson-Meyer:Descent} by
\begin{equation}
  \label{definitionofcoarsektheory}
  \KX_H^*(X,D) \defeq \Ktop_*\bigl(H, C_0(\Rips(X),D)\bigr).
\end{equation}
As observed in~\cite{Emerson-Meyer:Descent}, we have \(\Ktop_*\bigl(H,
C_0(\Rips(X),D)\bigr) \cong \K_*(C_0(\Rips(X),D)\cross H)\) because~\(H\) acts
properly on \(\Rips(X)\).

For most of our applications, \(X\) will be equivariantly uniformly
contractible for all compact subgroups \(K\subseteq H\), that is, the natural
embedding \(X\to\Rips(X)\) is a \(K\)\nb-equivariant coarse homotopy
equivalence.  In such cases, we simply have
\begin{equation}
  \label{eq:simplify_coarse_K-theory}%
  \KX_H^*(X,D) \cong \Ktop_*\bigl(H,C_0(X,D)\bigr).
\end{equation}
In particular, this applies if~\(X\) is an \(H\)\nb-compact universal proper
\(H\)\nb-space (again, recall that the coarse structure is determined by
requiring~\(H\) to act isometrically).

The \emph{\(H\)\nb-equivariant coarse co-assembly map for~\(X\) with
  coefficients in~\(D\)} is a certain map
\[
\mu^*\colon \Ktop_{*+1}\bigl(H,\bvr_H(X,D)\bigr) \to \KX^*_H(X,D)
\]
defined in~\cite{Emerson-Meyer:Descent}.  In the special case where we
have~\eqref{eq:simplify_coarse_K-theory}, this is simply the boundary map for
the extension \(C_0(X, D\otimes\Comp_H)\into \vvr_H(X,D) \prto \bvr_H(X,D)\).
We are implicitly using the fact that the functor \(\Ktop_*(H,\blank)\) has
long exact sequences for arbitrary extensions of \(H\)-\(C^*\)\nb-algebras,
which is proved in~\cite{Emerson-Meyer:Descent} using the isomorphism
\[
\Ktop_*(H,B) \defeq \K_*\bigl((B\otimes \ADir)\rcross H\bigr)
\cong \K_*\bigl((B\maxotimes \ADir)\cross H\bigr)
\]
and exactness properties of maximal \(C^*\)\nb-tensor products and full
crossed products.

There is also an alternative picture of the co-assembly map as a
forget-control map, provided~\(X\) is uniformly contractible
(see~\cite{Emerson-Meyer:Descent}*{\textsection2.8}).  We have the following
equivariant version of this result:

\begin{proposition}
  \label{forgettingcontrol}%
  Let~\(G\) be a totally disconnected group with a \(G\)\nb-compact
  universal proper \(G\)\nb-space \(\EG\).  Then the \(G\)\nb-equivariant
  coarse co-assembly map for~\(G\) is equivalent to the map
  \[
  j_*\colon \Ktop_{*+1}\bigl(G,\bvr_G(\EG,D)\bigr)
  \to\Ktop_{*+1}\bigl(G,\bbred_G(\EG,D)\bigr)
  \]
  induced by the inclusion \(j\colon \bvr_G(\EG,D) \to\bbred_G(\EG,D)\).
\end{proposition}

The equivalence of the two maps means that there is a natural commuting diagram
\[\xymatrix{
  \Ktop_{*+1}\bigl(G,\bvr_G(\abs{G},D)\bigr) \ar[r]^{\mu^*} \ar@{<->}[d]^{\cong} &
  \KX_*^G(\abs{G},D) \ar@{<->}[d]^{\cong} \\
  \Ktop_{*+1}\bigl(G,\bvr_G(\EG,D)\bigr) \ar[r]^{j^*} &
  \Ktop_{*+1}\bigl(G,\bbred_G(\EG,D)\bigr).
}
\]
Recall that~\(j\) is induced by the inclusion \(\vvr_G(\EG,D)\to
\bbbarred_G(\EG,D)\), which exactly forgets the vanishing variation condition.
Hence~\(j_*\) is a forget-control map.

\begin{proof}
  We may replace~\(\abs{G}\) by \(\EG\) because \(\EG\) is coarsely equivalent
  to~\(\abs{G}\).  The coarse \(\K\)\nb-theory of \(\EG\) agrees with the
  usual \(\K\)\nb-theory of \(\EG\) (see~\cite{Emerson-Meyer:Descent}).  A slight
  elaboration of the proof of \cite{Emerson-Meyer:Descent}*{Lemma 15} shows that
  \[
  \K^H_*\bigl(\bbbarred_G(\EG,D)\bigr) \cong
  \KK^H_*\bigl(\C,\bbbarred_G(\EG,D)\bigr)
  \]
  vanishes for all compact subgroups~\(H\) of~\(G\).  This yields
  \(\Ktop_*\bigl(G,\bbbarred_G(\EG,D)\bigr)=0\) by a result of~\cite{Chabert-Echterhoff-Oyono:Going}.
  Now the assertion follows from the Five Lemma and the naturality of the
  \(\K\)\nb-theory long exact sequence for~\eqref{eq:forget-control} as
  in~\cite{Emerson-Meyer:Descent}.
\end{proof}

\section{Classes in Kasparov theory from the stable Higson corona}

In this section, we show how to construct classes in equivariant
\(\KK\)-theory from the \(\K\)\nb-theory of the stable Higson corona.  The
following lemma is our main technical device:

\begin{lemma}
  \label{lem:coarse_gives_classes}%
  Let \(G\) and~\(H\) be locally compact groups and let~\(X\) be a coarse
  space equipped with commuting actions of \(G\) and~\(H\).  Suppose
  that~\(G\) acts by translations and that~\(H\) acts properly and by
  isometries.  Let \(A\) and~\(D\) be \(H\)\nb-\(C^*\)-algebras,
  equipped with the trivial \(G\)\nb-action.  We abbreviate
  \[
  B_X\defeq C_0(X,D\otimes\Comp_H\maxotimes A)\cross H,
  \qquad
  E_X\defeq (\vvr_H(X,D) \maxotimes A)\cross H
  \]
  and similarly for \(\Rips(X)\) instead of~\(X\).  There are extensions
  \(B_X\into E_X\prto E_X/B_X\) and \(B_{\Rips(X)}\into E_{\Rips(X)}\prto
  E_{\Rips(X)}/B_{\Rips(X)}\) with
  \[
  E_{\Rips(X)}/B_{\Rips(X)} \cong E_X/B_X
  \cong (\bvr_H(X,D) \maxotimes A)\cross H,
  \]
  and a natural commuting diagram
  \begin{equation}
    \label{eq:technical_diagram}
    \begin{gathered}
      \xymatrix{
        \K_{*+1}(E_X/B_X)
        \ar[r]^-{\partial} \ar[d]^{\psi} &
        \K_*(B_{\Rips(X)}) \ar[d]^{\phi} \\
        \KK^G_*(\C,B_X) \ar[r]^-{p_\EG^*} &
        \RKK^G_*(\EG;\C,B_X).
      }
    \end{gathered}
  \end{equation}
\end{lemma}

\begin{proof}
  The quotients \(E_X/B_X\) and \(E_{\Rips(X)}/B_{\Rips(X)}\) are as asserted
  and agree because \(X\to\Rips(X)\) is a coarse equivalence and because
  maximal tensor products and full crossed products are exact functors in
  complete generality, unlike spatial tensor products and reduced crossed
  products.  We let~\(\partial\) be the \(\K\)\nb-theory boundary map for
  the extension \(B_{\Rips(X)}\into E_{\Rips(X)}\prto E_X/B_X\).

  Since we have a natural map \(\bvr_H(X,D)\maxotimes A\to
  \bvr_H(X,D\maxotimes A)\), we may replace the pair \((D,A)\) by
  \((D\maxotimes A,\C)\) and omit~\(A\) if convenient.  Stabilising~\(D\)
  by~\(\Comp_H\), we can further eliminate the stabilisations.

  First we lift the \(\K\)\nb-theory boundary map for the extension
  \(B_X\into E_X\prto E_X/B_X\) to a map \(\psi\colon \K_{*+1}(E_X/B_X)\to
  \KK^G_*(\C,B_X)\).  The \(G\)\nb-equivariance of the resulting Kasparov
  cycles follows from the assumption that~\(G\) acts on~\(X\) by translations.

  We have to distinguish between the cases \(*=0\) and \(*=1\).  We only write
  down the construction for \(*=0\).  Since the algebra \(E_X/B_X\) is
  matrix-stable, \(\K_1(E_X/B_X)\) is the homotopy group of unitaries in
  \(E_X/B_X\) without further stabilisation.  A cycle for \(\KK^G_0(\C,B_X)\)
  is given by two \(G\)\nb-equivariant Hilbert modules~\(\Hilm_\pm\)
  over~\(B_X\) and a \(G\)\nb-continuous adjointable operator
  \(F\colon\Hilm_+\to\Hilm_-\) for which \(1-FF^*\), \(1-F^*F\) and \(gF-F\)
  for \(g\in G\) are compact; we take \(\Hilm_\pm = B_X\) and let \(F\in
  E_X\subseteq \Mult(B_X)\) be a lifting for a unitary \(u\in E_X/B_X\).
  Since~\(G\) acts on~\(X\) by translations, the induced action on
  \(\vvr(X,D)\) and hence on \(E_X/B_X\) is trivial.  Hence~\(u\) is a
  \(G\)\nb-invariant unitary in \(E_X/B_X\).  For the lifting~\(F\), this
  means that
  \[
  1-FF^*,\ 1-F^*F,\ gF-F\in B_X.
  \]
  Hence~\(F\) defines a cycle for \(\KK_0^G(\C,B_X)\).  We get a well-defined
  map \([u]\mapsto [F]\) from \(\K_1(E_X/B_X)\) to \(\KK^G_0(\C,B_X)\) because
  homotopic unitaries yield operator homotopic Kasparov cycles.

  Next we have to factor the map \(p_\EG^*\circ\psi\)
  in~\eqref{eq:technical_diagram} through \(\K_0(B_{\Rips(X)})\).  The main
  ingredient is a certain continuous map \(\bar{c}\colon \EG\times X\to\Rips(X)\).
  We use the same description of \(\Rips(X)\) as in~\cite{Emerson-Meyer:Descent} as
  the space of positive measures on~\(X\) with \(1/2 < \mu(X)\le1\); this is a
  \(\sigma\)\nb-coarse space in a natural way, we write it as
  \(\Rips(X)=\bigcup \Subrips_d(X)\).

  There is a function \(c\colon \EG\to\R_+\) for which
  \(\int_\EG c(\mu g)\,\textup{d}g=1\) for all \(\mu\in\EG\)
  and \(\supp c\cap Y\) is compact for \(G\)\nb-compact
  \(Y\subseteq\EG\).  If \(\mu\in\EG\), \(x\in X\), then the
  condition
  \[
  \langle \bar{c}(\mu,x), \alpha\rangle \defeq
  \int_G c(\mu g) \alpha(g^{-1}x)\,\textup{d}g
  \]
  for \(\alpha\in C_0(X)\) defines a probability measure on~\(X\).  Since such
  measures are contained in \(\Rips(X)\), \(\bar{c}\) defines a map
  \(\bar{c}\colon \EG\times X\to\Rips(X)\).  This map is continuous and satisfies
  \(\bar{c}(\mu g,g^{-1}xh)=\bar{c}(\mu,x)h\) for all \(g\in G\),
  \(\mu\in\EG\), \(x\in X\), \(h\in H\).

  For a \(C^*\)\nb-algebra~\(Z\), let \(C(\EG,Z)\) be the
  \(\sigma\)\nb-\(C^*\)-algebra of all continuous functions \(f\colon \EG\to
  Z\) without any growth restriction.  Thus \(C(\EG,Z)=\varprojlim C(K,Z)\),
  where~\(K\) runs through the directed set of compact subsets of~\(\EG\).

  We claim that \((\bar{c}^* f)(\mu)(x) \defeq f\bigl(\bar{c}(\mu,x)\bigr)\)
  for \(f\in C_0(\Rips(X),D)\) defines a continuous \(*\)\nb-homomorphism
  \[
  \bar{c}^*\colon C_0(\Rips(X),D) \to C\bigl(\EG,C_0(X,D)\bigr).
  \]
  If \(K\subseteq\EG\) is compact, then there is a compact subset \(L\subseteq
  G\) such that \(c(\mu\cdot g)=0\) for \(\mu\in K\) and \(g\notin L\).  Hence
  \(\bar{c}(\mu,x)\) is supported in \(L^{-1}x\) for \(\mu\in K\).
  Since~\(G\) acts on~\(X\) by translations, such measures are contained in a
  filtration level \(\Subrips_d(X)\).  Hence \(\bar{c}^*(f)\) restricts to a
  \(C_0\)-function \(K\times X\to D\) for all \(f\in C_0(\Rips(X),D)\).  This proves
  the claim.  Since~\(\bar{c}\) is \(H\)\nb-equivariant and
  \(G\)\nb-invariant, we get an induced map
  \[
  B_{\Rips(X)} = C_0(\Rips(X),D)\cross H \to
  (C\bigl(\EG,C_0(X,D)\bigr)\cross H)^G = C(\EG,B_X)^G,
  \]
  where \(Z^G\subseteq Z\) denotes the subalgebra of \(G\)\nb-invariant
  elements.  We obtain an induced \(*\)\nb-homomorphism between the stable
  multiplier algebras as well.

  An element of \(\K_0\bigl(B_{\Rips(X)})\) is represented by a self-adjoint
  bounded multiplier \(F\in\Mult(B_{\Rips(X)}\otimes\Comp)\) such that
  \(1-FF^*\) and \(1-F^*F\) belong to \(B_{\Rips(X)}\otimes\Comp\).  Now
  \(\tilde{F}\defeq \bar{c}^*(F)\) is a \(G\)\nb-invariant bounded
  multiplier of \(C(\EG,B_X\otimes\Comp)\) and hence a \(G\)\nb-invariant
  multiplier of \(C_0(\EG,B_X\otimes\Comp)\), such that \(\alpha\cdot
  (1-\tilde{F}\tilde{F}^*)\) and \(\alpha\cdot (1-\tilde{F}^*\tilde{F})\)
  belong to \(C_0(\EG,B_X\otimes\Comp)\) for all \(\alpha\in C_0(\EG)\).  This
  says exactly that~\(\tilde{F}\) is a cycle for \(\RKK^G_0(\EG;\C,B_X)\).
  This construction provides the natural map
  \[
  \phi\colon \K_0(B_{\Rips(X)}) \to \RKK^G_0(\EG;\C,B_X).
  \]

  Finally, a routine computation, which we omit, shows that the two images of
  a unitary \(u\in E_X/B_X\) differ by a compact perturbation.  Hence the
  diagram~\eqref{eq:technical_diagram} commutes.
\end{proof}

We are mainly interested in the case where \(A\) is the source~\(\ADir\) of a
Dirac morphism for~\(H\).  Then \(\K_{*+1}(E_X/B_X) =
\Ktop_*\bigl(H,\bvr_H(X,D)\bigr)\), and the top row
in~\eqref{eq:technical_diagram} is the \(H\)\nb-equivariant coarse
co-assembly map for~\(X\) with coefficients~\(D\).  Since we assume~\(H\) to
act properly on~\(X\), we have a \(\KK^G\)-equivalence \(B_X\sim
C_0(X,D)\cross H\), and similarly for \(\Rips(X)\).  Hence we now get a
commuting square
\begin{equation}
  \label{eq:notso_technical_diagram}
  \begin{gathered}
    \xymatrix@C=4em{
      \Ktop_{*+1}\bigl(H,\bvr_H(X,D)\bigr) \ar[r]^-{\partial}
      \ar[d]^{\psi} &
      \KX^*_H(X,D) \ar[d]^{\phi} \\
      \KK^G_*(\C,C_0(X,D)\cross H) \ar[r]^-{p_\EG^*} &
      \RKK^G_*(\EG;\C,C_0(X,D)\cross H).
    }
  \end{gathered}
\end{equation}
If, in addition, \(D=\C\) and the action of~\(H\) on~\(X\) is free, then we
can further simplify this to
\begin{equation}
  \label{eq:even_less_technical_diagram}
  \begin{gathered}
    \xymatrix@C=4em{
      \Ktop_{*+1}\bigl(H,\bvr_H(X)\bigr) \ar[r]^-{\partial} \ar[d]^{\psi} &
      \KX^*_H(X) \ar[d]^{\phi} \\
      \KK^G_*\bigl(\C,C_0(X/H)\bigr) \ar[r]^-{p_\EG^*} &
      \RKK^G_*\bigl(\EG;\C,C_0(X/H)\bigr).
    }
  \end{gathered}
\end{equation}

We may also specialise the space~\(X\) to~\(\abs{G}\), with~\(G\) acting by
multiplication on the left, and with \(H\subseteq G\) a compact subgroup
acting on~\(\abs{G}\) by right multiplication.  This is the special case
of~\eqref{eq:technical_diagram} that is used in~\cite{Emerson-Meyer:Descent}.  The
following applications will require other choices of~\(X\).

\subsection{Applications to Lipschitz classes}
\label{sec:apply_main_lemma}

Now we use Lemma~\ref{lem:coarse_gives_classes} to construct interesting
elements in \(\KK^G_*\bigl(\C,C_0(X)\bigr)\) for a \(G\)\nb-space~\(X\).
This is related to the method of Lipschitz maps developed by Connes, Gromov
and Moscovici in~\cite{Connes-Gromov-Moscovici:Lipschitz}.

\subsubsection{Pulled-back coarse structures}
\label{sec:pull-back_coarse}

Let~\(X\) be a \(G\)\nb-space, let~\(Y\) be a coarse space and let
\(\alpha\colon X\to Y\) be a proper continuous map.  We pull back the coarse
structure on~\(Y\) to a coarse structure on~\(X\), letting \(E\subseteq X\times X\)
be an entourage if and only if \(\alpha_*(E)\subseteq Y\times Y\) is one.
Since~\(\alpha\) is proper and continuous, this coarse structure is compatible
with the topology on~\(X\).  For this coarse structure, \(G\) \emph{acts by
  translations} if and only if~\(\alpha\) satisfies the following
\emph{displacement condition} used in~\cite{Connes-Gromov-Moscovici:Lipschitz}: for any compact subset
\(K\subseteq G\), the set
\[
\bigl\{\bigl(\alpha(gx),\alpha(x)\bigr)\in Y\times Y \bigm|
x\in X,\ g\in K\bigr\}
\]
is an entourage of~\(Y\).  The map~\(\alpha\) becomes a coarse map.  Hence we
obtain a commuting diagram
\[\xymatrix{
  \K_{*+1}\bigl(\bvr(Y)\bigr) \ar[r]^{\alpha^*} \ar[d]^{\partial^Y} &
  \K_{*+1}\bigl(\bvr(X)\bigr) \ar[r]^-{\psi} \ar[d]^{\partial^X} &
  \KK^G_*\bigl(\C,C_0(X)\bigr) \ar[d]^{p_\EG^*} \\
  \KX^*(Y) \ar[r]^{\alpha^*} &
  \KX^*(X) \ar[r]^-{\phi} & \RKK^G_*\bigl(\EG;\C,C_0(X)\bigr).
}
\]
with \(\psi\) and~\(\phi\) as in Lemma~\ref{lem:coarse_gives_classes}.

The constructions of~\cite{Connes-Gromov-Moscovici:Lipschitz}*{\textsection I.10} only use \(Y=\R^N\) with
the Euclidean coarse structure.  The coarse co\nb-assembly map is an
isomorphism for~\(\R^N\) because~\(\R^N\) is scalable.  Moreover, \(\R^N\) is
uniformly contractible and has bounded geometry.  Hence we obtain canonical
isomorphisms
\[
\K_{*+1}\bigl(\bvr(\R^N)\bigr) \cong \KX^*(\R^N) \cong \K^*(\R^N).
\]
In particular, \(\K_{*+1}\bigl(\bvr(\R^N)\bigr)\cong\Z\) with generator
\([\partial\R^N]\) in \(\K_{1-N}\bigl(\bvr(\R^N)\bigr)\).  This class is
nothing but the usual dual-Dirac morphism for the locally compact
group~\(\R^N\).  As a result, any map \(\alpha\colon X\to\R^N\) that satisfies
the displacement condition above induces
\[
[\alpha]\defeq \psi\bigl(\alpha^*[\partial\R^N]\bigr)
\in\KK^G_{-N}\bigl(\C,C_0(X)\bigr).
\]
The commutative diagram~\eqref{eq:technical_diagram} computes
\(p_\EG^*[\alpha]\in \RKK^G_{-N}\bigl(\EG;\C,C_0(X)\bigr)\) in purely
topological terms.

\subsubsection{Principal bundles over coarse spaces}
\label{sec:bundles_over_coarse}

As in~\cite{Connes-Gromov-Moscovici:Lipschitz}, we may replace a fixed map \(X\to\R^N\) by a section of a
vector bundle over~\(X\).  But we need this bundle to have a
\(G\)\nb-equivariant spin structure.  To encode this, we consider a
\(G\)\nb-equivariant \(\Spin(N)\)\nb-principal bundle \(\pi\colon E\to B\)
together with actions of~\(G\) on \(E\) and~\(B\) such that~\(\pi\) is
\(G\)\nb-equivariant and the action on~\(E\) commutes with the action of
\(H\defeq \Spin(N)\).  Let \(T\defeq E\times_{\Spin(N)} \R^N\) be the associated
vector bundle over~\(B\).  It carries a \(G\)\nb-invariant Euclidean metric
and spin structure.  As is well-known, sections \(\alpha\colon B\to T\)
correspond bijectively to \(\Spin(N)\)-equivariant maps \(\alpha'\colon
E\to\R^N\); here a section~\(\alpha\) corresponds to the map \(\alpha'\colon
E\to\R^N\) that sends \(y\in E\) to the coordinates of \(\alpha\pi(y)\) in the
orthogonal frame described by~\(y\).  Since the group \(\Spin(N)\) is compact,
the map~\(\alpha'\) is proper if and only if \(b\mapsto \norm{\alpha(b)}\) is
a proper function on~\(B\).

As in \textsection\ref{sec:pull-back_coarse}, a \(\Spin(N)\)\nb-equivariant
proper continuous map \(\alpha'\colon E\to Y\) for a coarse space~\(Y\) allows
us to pull back the coarse structure of~\(Y\) to~\(E\); then \(\Spin(N)\) acts
by isometries.  The group~\(G\) acts by translations if and only
if~\(\alpha'\) satisfies the displacement condition
from~\textsection\ref{sec:pull-back_coarse}.  If \(Y=\R^N\), we can rewrite
this in terms of \(\alpha\colon B\to T\): we need
\[
\sup \bigl\{ \norm{g\alpha(g^{-1}b)-\alpha(b)} \bigm|
b\in B,\ g\in K \bigr\}
\]
to be bounded for all compact subsets \(K\subseteq G\).

If the displacement condition holds, then we are in the
situation of Lemma~\ref{lem:coarse_gives_classes} with
\(H=\Spin(N)\) and \(X=E\).  Since~\(H\) acts freely on~\(E\),
\(C_0(E)\cross H\) is \(G\)\nb-equivariantly Morita--Rieffel
equivalent to \(C_0(B)\).  We obtain canonical maps
\begin{multline*}
  \K^{\Spin(N)}_{*+1}\bigl(\bvr_{\Spin(N)}(\R^N)\bigr)
  \xrightarrow{(\alpha')^*}
  \K^{\Spin(N)}_{*+1}\bigl(\bvr_{\Spin(N)}(E)\bigr) \\
  \xrightarrow{\psi} \KK^G_*\bigl(\C,C_0(E)\cross \Spin(N)\bigr)
  \cong \KK^G_*\bigl(\C,C_0(B)\bigr).
\end{multline*}

The \(\Spin(N)\)\nb-equivariant coarse co-assembly map for~\(\R^N\) is an
isomorphism by~\cite{Emerson-Meyer:Descent} because the group \(\R^N\cross\Spin(N)\)
has a dual-Dirac morphism.  Using also the uniform contractibility of~\(\R^N\)
and \(\Spin(N)\)-equivariant Bott periodicity, we get
\[
\K^{\Spin(N)}_{*+1}\bigl(\bvr_{\Spin(N)}(\R^N)\bigr)
\cong \KX_{\Spin(N)}^*(\R^N)
\cong \K_{\Spin(N)}^*(\R^N)
\cong \K_{\Spin(N)}^{*+N}(\point).
\]
The class of the trivial representation in \(\Rep(\Spin N) \cong
\K^{\Spin(N)}_0(\C)\) is mapped to the usual dual-Dirac morphism \([\partial
\R^N]\in \K^{\Spin(N)}_{1-N}\bigl(\bvr_{\Spin(N)}(\R^N)\bigr)\) for~\(\R^N\).
As a result, any proper section \(\alpha\colon B\to T\) satisfying the
displacement condition induces
\[
[\alpha]\defeq \psi\bigl(\alpha^*[\partial\R^N]\bigr) \in
\KK^G_{-N}\bigl(\C,C_0(B)\bigr).
\]
Again, the commutative diagram~\eqref{eq:technical_diagram} computes
\(p_\EG^*[\alpha]\in \RKK^G_{-N}\bigl(\EG;\C,C_0(X)\bigr)\) in purely
topological terms.

\subsubsection{Coarse structures on jet bundles}
\label{sec:coarse_jet}

Let~\(M\) be an oriented compact manifold and let \(\Diff^+(M)\) be the
infinite-dimensional Lie group of orientation-preserving diffeomorphisms
of~\(M\).  Let~\(G\) be a locally compact group that acts on~\(M\) by a
continuous group homomorphism \(G\to\Diff^+(M)\).  The \emph{Gelfand--Fuchs
  cohomology} of~\(M\) is part of the group cohomology of \(\Diff^+(M)\) and
by functoriality maps to the group cohomology of~\(G\).  It is shown
in~\cite{Connes-Gromov-Moscovici:Lipschitz} that the range of Gelfand--Fuchs cohomology in the cohomology
of~\(G\) yields homotopy-invariant higher signatures.

This argument has two parts; one is geometric and concerns the construction of
a class in \(\KK^G_*\bigl(\C,C_0(X)\bigr)\) for a suitable space~\(X\); the
other uses cyclic homology to construct linear functionals on
\(\K_*(C_0(X)\rcross G)\) associated to Gelfand--Fuchs cohomology classes.  We
can simplify the first step; the second has nothing to do with coarse
geometry.

Let \(\pi^k\colon \Jet^k_+(M)\to M\) be the \emph{oriented \(k\)\nb-jet
  bundle} over~\(M\).  That is, a point in \(\Jet^k_+(M)\) is the \(k\)th
order Taylor series at~\(0\) of an orientation-preserving diffeomorphism from
a neighbourhood of \(0\in\R^n\) into~\(M\).  This is a principal
\(H\)\nb-bundle over~\(M\), where~\(H\) is a connected Lie group whose Lie
algebra~\(\mathfrak{h}\) is the space of polynomial maps \(p\colon
\R^n\to\R^n\) of order~\(k\) with \(p(0)=0\), with an appropriate Lie algebra
structure.  The maximal compact subgroup \(K\subseteq H\) is isomorphic to
\(\textup{SO}(n)\), acting by isometries on~\(\R^n\).  It acts
on~\(\mathfrak{h}\) by conjugation.

Since our construction is natural, the action of~\(G\) on~\(M\) lifts to an
action on \(\Jet^k_+(M)\) that commutes with the \(H\)\nb-action.  We
let~\(H\) act on the right and~\(G\) on the left.  Define \(X_k\defeq
\Jet^k_+(M)/K\).  This is the bundle space of a fibration over~\(M\) with
fibres \(H/K\).  Gelfand--Fuchs cohomology can be computed using a chain
complex of \(\Diff^+(M)\)-invariant differential forms on~\(X_k\) for
\(k\to\infty\).  Using this description, Connes, Gromov, and Moscovici
associate to a Gelfand--Fuchs cohomology class a functional
\(\K_*(C_0(X_k)\rcross G)\to\C\) for sufficiently high~\(k\) in~\cite{Connes-Gromov-Moscovici:Lipschitz}.

Since \(\Jet^k_+(M)/H\cong M\) is compact, there is a unique coarse structure
on \(\Jet^k_+(M)\) for which~\(H\) acts isometrically (see
\textsection\ref{sec:acts_coarse}).  With this coarse structure,
\(\Jet^k_+(M)\) is coarsely equivalent to~\(H\).  The compactness of
\(\Jet^k_+(M)/H\cong M\) also implies easily that~\(G\) acts by translations.
We have a Morita--Rieffel equivalence \(C_0(X_k)\sim C_0(\Jet^k_+ M)\cross K\)
because~\(K\) acts freely on \(\Jet^k_+(M)\).  We want to study the map
\[
\K_{*+1}(\bvr_K(\Jet^k_+ M)\cross K) \xrightarrow{\psi}
\KK^G_*(\C,C_0(\Jet^k_+ M)\cross K) \cong
\KK^G_*\bigl(\C,C_0(X_k)\bigr)
\]
produced by Lemma~\ref{lem:coarse_gives_classes}.

Since~\(H\) is almost connected, it has a dual-Dirac morphism
by~\cite{Kasparov:Novikov}; hence the \(K\)\nb-equivariant coarse co-assembly map
for~\(H\) is an isomorphism by the main result of~\cite{Emerson-Meyer:Descent}.
Moreover, \(H/K\) is a model for \(\EG\) by~\cite{Abels:Slices}.  We get
\[
\K_{*+1}(\bvr_K(\Jet^k_+ M)\cross K)
\cong \K_{*+1}(\bvr_K(\abs{H})\cross K)
\cong \KX^*_K(\abs{H})
\cong \K^*_K(H/K).
\]
Let \(\mathfrak{h}\) and~\(\mathfrak{k}\) be the Lie algebras of \(H\)
and~\(K\).  There is a \(K\)\nb-equivariant homeomorphism
\(\mathfrak{h}/\mathfrak{k}\cong H/K\), where~\(K\) acts on
\(\mathfrak{h}/\mathfrak{k}\) by conjugation.  Now we need to know whether
there is a \(K\)\nb-equivariant spin structure on
\(\mathfrak{h}/\mathfrak{k}\).  One can check that this is the case if
\(k\equiv 0,1 \bmod 4\).  Since we can choose~\(k\) as large as we like, we
can always assume that this is the case.  The spin structure allows us to use
Bott periodicity to identify \(\K^*_K(H/K) \cong \K_{*-N}^H(\C)\), which is
the representation ring of~\(K\) in degree~\(-N\), where \(N=\dim
\mathfrak{h}/\mathfrak{k}\).  Using our construction, the trivial
representation of~\(K\) yields a canonical element in
\(\KK^G_{-N}\bigl(\C,C_0(X_k)\bigr)\).

This construction is much shorter than the corresponding one in~\cite{Connes-Gromov-Moscovici:Lipschitz}
because we use Kasparov's result about dual-Dirac morphisms for almost
connected groups.  Much of the corresponding argument in~\cite{Connes-Gromov-Moscovici:Lipschitz} is
concerned with proving a variation on this result of Kasparov.

\section{Computation of \texorpdfstring{$\KK^G_*(\C,\ADir)$}{KKG(C,P)}}
\label{sec:compute_KKGCP}

So far, we have merely used the diagram~\eqref{eq:technical_diagram} to
construct certain elements in \(\KK^G_*(\C,B)\).  Now we show that this
construction yields an isomorphic description of \(\KK^G_*(\C,\ADir)\).  This
assertion requires~\(G\) to be a totally disconnected group with a
\(G\)\nb-compact universal proper \(G\)\nb-space.  We assume this
throughout this section.

\begin{lemma}
  \label{hatsoff}%
  In the situation of
  Lemma~\textup{\ref{lem:coarse_gives_classes}}, suppose that
  \(X = \abs{G}\) with~\(G\) acting by left translations and
  that \(H\subseteq G\) is a compact subgroup acting on~\(X\)
  by right translations; here~\(\abs{G}\) carries the coarse
  structure of Example~\textup{\ref{groupsascoarsespaces}}.
  Then the maps \(\psi\) and~\(\phi\) are isomorphisms.
\end{lemma}

\begin{proof}
  We reduce this assertion to results of~\cite{Emerson-Meyer:Descent}.  The
  \(C^*\)\nb-algebras \(\bvr_H(\abs{G},D)\cross H\) and \(\bvr_H(\abs{G},
  D)^H\) are strongly Morita equivalent, whence have isomorphic
  \(\K\)\nb-theory.  It is shown in~\cite{Emerson-Meyer:Descent} that
  \begin{equation}
    \label{car}
    \K_{*+1}\bigl(\bvr_H(\abs{G},D)^H\bigr) \cong \KK^G_*(\C,\Ind_H^G D).
  \end{equation}
  Finally, \(\Ind_H^G(D) = C_0(G,D)^H\) is \(G\)\nb-equivariantly
  Morita--Rieffel equivalent to \(C_0(G,D)\cross H\).  Hence we get
  \begin{multline}
    \K_{*+1}(\bvr_H(\abs{G},D)\cross H)
    \cong \K_{*+1}\bigl(\bvr_H(\abs{G},D)^H\bigr)
    \cong \KK^G_*\bigl(\C, \Ind_H^G(D)\bigr)
    \\ \cong\KK^G_*\bigl(\C, C_0(G,D)\cross H\bigr).
  \end{multline}
  It is a routine exercise to verify that this composition agrees with the
  map~\(\psi\) in~\eqref{eq:technical_diagram}.  Similar considerations apply
  to the map~\(\phi\).
\end{proof}

We now set \(X=\abs{G}\), and let \(G=H\).  The actions of~\(G\)
on~\(\abs{G}\) on the left and right are by translations and isometries,
respectively.  Lemma~\ref{lem:coarse_gives_classes} yields a map
\begin{equation}
  \label{map}
  \Psi_*\colon
  \K_{*+1}\bigl((\bvr(\abs{G},D)\maxotimes A)\cross G\bigr) \to
  \KK^G_*\bigl(\C,C_0(\abs{G},D\maxotimes A)\cross G\bigr).
\end{equation}
for all~\(A,D\), where we use the \(G\)-equivariant Morita--Rieffel equivalence
between \(C_0(\abs{G},D)\cross G\) and~\(D\).  It fits into a commuting
diagram
\[\xymatrix@C=4em{
  \K_{*+1}\bigl((\bvr_G(\abs{G},D)\maxotimes A)\cross G\bigr)
  \ar[r]^-{\partial}
  \ar[d]^{\Psi^{D,A}_*} &
  \K_*(C_0(\EG,D\maxotimes A)\cross G) \ar[d]_{\cong} \\
  \KK_*^G(\C,D\maxotimes A) \ar[r]^-{p_\EG^*} &
  \RKK_*^G(\EG;\C,D\maxotimes A).
}
\]

\begin{lemma}
  \label{lem:Psi_iso}%
  The class of \(G\)-\(C^*\)-algebras~\(A\) for which \(\Psi_*^{D,A}\) is an
  isomorphism for all~\(D\) is triangulated and thick and contains all
  \(G\)-\(C^*\)\nb-algebras of the form \(C_0(G/H)\) for compact open
  subgroups~\(H\) of~\(G\).
\end{lemma}

\begin{proof}
  The fact that this category of algebras is triangulated and thick means that
  it is closed under suspensions, extensions, and direct summands.  These
  formal properties are easy to check.

  Since \(H\subseteq G\) is open, there is no difference between
  \(H\)\nb-continuity and \(G\)\nb-\hspace{0pt}continuity.  Hence
  \begin{align*}
    \bigl(\bvr_G(\abs{G},D)\maxotimes C_0(G/H)\bigr)\cross G
    &\cong \bigl(\bvr_H(\abs{G},D)\maxotimes C_0(G/H)\bigr)\cross G
    \\ &\sim \bigl(\bvr_H(\abs{G},D) \cross H,
  \end{align*}
  where~\(\sim\) means Morita--Rieffel equivalence.  Similar simplifications
  can be made in other corners of the square.  Hence the diagram for
  \(A=C_0(G/H)\) and~\(G\) acting on the right is equivalent to a
  corresponding diagram for trivial~\(A\) and~\(H\) acting on the right.  The
  latter case is contained in Lemma~\ref{hatsoff}.
\end{proof}

\begin{theorem}
  \label{the:Kreg}
  Let~\(G\) be an almost totally disconnected group with \(G\)\nb-compact
  \(\EG\).  Then for every \(B\in\KK^G\), the map
  \[
  \Psi_*\colon \Ktop_{*+1}\bigl(G,\bvr(\abs{G},B)\bigr) \to
  \KK_*^G(\C,B\otimes\ADir)
  \]
  is an isomorphism and the diagram
  \[
  \xymatrix{
    \Ktop_{*+1}\bigl(G,\bvr(\abs{G},B)\bigr)
    \ar[r]^-{\mu^*} \ar[d]_{\cong}^{\Psi_*} &
    \KX^*_G(\abs{G},B) \ar[d]_{\cong} \\
    \KK_*^G(\C,B\otimes\ADir) \ar[r]^-{p_\EG^*} &
    \RKK_*^G(\EG;\C,B\otimes\ADir)
  }
  \]
  commutes.  In particular, \(\Ktop_{*+1}\bigl(G,\bvr(\abs{G})\bigr)\) is
  naturally isomorphic to \(\KK^G_*(\C,\ADir)\).
\end{theorem}

\begin{proof}
  It is shown in~\cite{Emerson-Meyer:Descent} that for such groups~\(G\), the
  algebra~\(\ADir\) belongs to the thick triangulated subcategory of \(\KK^G\)
  that is generated by \(C_0(G/H)\) for compact subgroups~\(H\) of~\(G\).
  Hence the assertion follows from Lemma~\ref{lem:Psi_iso} and our definition
  of \(\Ktop\).
\end{proof}

\begin{corollary}
  \label{cor:coassembly-equivalent}%
  Let \(\Dirac\in\KK^G(\ADir,\C)\) be a Dirac morphism for~\(G\).  Then the
  following diagram commutes
  \[
  \xymatrix{
    \Ktop_{*+1}\bigl(G,\bvr_G(\abs{G})\bigr)
    \ar@/^3em/[rr]^-{\mu^*}
    \ar[d]^{\Psi_*}_{\cong} \ar[r]^-{\partial} &
    \K_*(C_0(\EG,\ADir)\cross G) \ar[d]_{\cong}
    \ar[r]^{\Dirac_*}_{\cong} &
    \K_*(C_0(\EG)\cross G) \ar[d]_{\cong} \\
    \KK_*^G(\C,\ADir) \ar[r]^-{p_\EG^*} \ar[dr]_-{\Dirac_*} &
    \RKK_*^G(\EG;\C,\ADir) \ar[r]^{\Dirac_*}_{\cong} &
    \RKK_*^G(\EG;\C,\C) \\
    & \KK^G_*(\C,\C) \ar[ur]_-{p_\EG^*},
  }
  \]
  where~\(\Psi_*\) is as in Theorem~\textup{\ref{the:Kreg}}
  and~\(\mu^*\) is the \(G\)\nb-equivariant coarse co-assembly
  map for~\(\abs{G}\), and the indicated maps are isomorphisms.
\end{corollary}

\begin{proof}
  This follows from Theorem~\ref{the:Kreg} and the general properties of the
  Dirac morphism discussed in~\textsection\ref{sec:Dirac_and_BC}.
\end{proof}

\begin{definition}
  We call \(a\in \RKK^G_*(\EG;\C,\C)\)
  \begin{enumerate}[label=\textup{(\alph*)}]
  \item a \emph{boundary class} if it lies in the range of
    \[
    \mu^*\colon \Ktop_{*+1}\bigl(G,\bvr(\abs{G})\bigr)\to \RKK^G_*(\EG;\C,\C);
    \]

  \item \emph{properly factorisable} if \(a = p_\EG^*(b\otimes_A c)\) for some
    proper \(G\)-\(C^*\)-algebra~\(A\) and some \(b\in\KK^G_*(\C,A)\),
    \(c\in\KK^G_*(A,\C)\);

  \item \emph{proper Lipschitz} if \(a = p_\EG^*(b\otimes_{C_0(X)}c)\), where
    \(b\in\KK^G_*\bigl(\C,C_0(X))\) is constructed as in
    \textsection\ref{sec:pull-back_coarse} and
    \textsection\ref{sec:bundles_over_coarse} and \(c\in \KK^G_*(C_0(X),\C)\)
    is arbitrary.
 \end{enumerate}
\end{definition}

\begin{proposition}
  Let~\(G\) be a totally disconnected group with \(G\)\nb-compact \(\EG\).
  \begin{enumerate}[label=\textup{(\alph*)}]
  \item A class \(a\in \RKK^G_*(\EG;\C,\C)\) is properly factorisable if and
    only if it is a boundary class.

  \item Proper Lipschitz classes are boundary classes.

  \item The boundary classes form an ideal in the ring
    \(\RKK^G_*(\EG;\C,\C)\).

  \item The class \(1_\EG\) is a boundary class if and only if~\(G\) has a
    dual-Dirac morphism; in this case, the \(G\)\nb-equivariant coarse
    co-assembly map~\(\mu^*\) is an isomorphism.

  \item Any boundary class lies in the range of
    \[
    p_\EG^*\colon \KK^G_*(\C,\C)\to\RKK^G_*(\EG;\C,\C)
    \]
    and hence yields homotopy invariants for manifolds.
  \end{enumerate}
\end{proposition}

\begin{proof}
  By Corollary~\ref{cor:coassembly-equivalent}, the equivariant coarse
  co-assembly map
  \[
  \mu^*\colon
  \Ktop_{*+1}\bigl(G,\bvr(\abs{G})\bigr) \to
  \K_*\bigl(C_0(\EG\cross G)\bigr) \cong \RKK^G_*(\EG;\C,\C)
  \]
  is equivalent to the map
  \[
  \KK^G_*(\C,\ADir) \xrightarrow{p_\EG^*}
  \RKK^G_*(\EG;\C,\ADir) \cong
  \RKK^G_*(\EG;\C,\C).
  \]
  If we combine this with the isomorphism \(\RKK^G_*(\EG;\C,\C) \cong
  \KK^G_*(\ADir,\ADir)\), the resulting map
  \[
  \KK^G_*(\C,\ADir) \to \KK^G_*(\ADir,\ADir)
  \]
  is simply the product (on the left) with \(\Dirac\in\KK^G(\ADir,\C)\); this
  map is known to be an isomorphism if and only if it is surjective, if and
  only if~\(1_\ADir\) is in its range, if and only if the
  \(H\)\nb-equivariant coarse co-assembly map
  \[
  \K_*(\bvr_H(\abs{G})\cross H) \to \KX^*_H(\abs{G})
  \]
  is an isomorphism for all compact subgroups~\(H\) of~\(G\)
  by~\cite{Emerson-Meyer:Descent}.

  For any \(G\)\nb-\(C^*\)-algebra~\(B\), the \(\Ztwo\)\nb-graded group
  \(\Ktop_*(G,B)\cong \K_*\bigl((B\otimes\ADir)\rcross G\bigr)\) is a graded
  module over the \(\Ztwo\)\nb-graded ring \(\KK^G_*(\ADir,\ADir)\cong
  \RKK^G_*(\EG;\C,\C)\) in a canonical way; the isomorphism between these two
  groups is a ring isomorphism because it is the composite of the two ring
  isomorphisms
  \[
  \KK^G_*(\ADir,\ADir)\xrightarrow[\cong]{p_\EG^*}
  \RKK^G_*(\EG;\ADir,\ADir)\xleftarrow[\cong]{\blank\otimes\ADir}
  \RKK^G_*(\EG;\C,\C).
  \]
  Hence we get module structures on \(\Ktop_*\bigl(G,\bvr(\abs{G})\bigr)\) and
  \[
  \Ktop_*\bigl(G,C_0(\EG,\Comp)\bigr) \cong
  \K_*(C_0(\EG)\cross G) \cong
  \RKK^G_*(\EG;\C,\C).
  \]
  The latter isomorphism is a module isomorphism; thus
  \(\Ktop_*\bigl(G,C_0(\EG,\Comp)\bigr)\) is a free module of rank~\(1\).  The
  equivariant co-assembly map is natural in the formal sense, so that it is an
  \(\RKK^G_*(\EG;\C,\C)\)-module homomorphism.  Hence its range is a
  submodule, that is, an ideal in \(\RKK^G_*(\EG;\C,\C)\) (since this ring is
  graded commutative, there is no difference between one- and two-sided graded
  ideals).  This also yields that~\(\mu^*\) is surjective if and only if it is
  bijective, if and only if the unit class \(1_\EG\) belongs to its range; we
  already know this from~\cite{Emerson-Meyer:Descent}.

  If~\(A\) is a proper \(G\)-\(C^*\)\nb-algebra, then
  \(\ID_A\otimes\Dirac\in\KK^G(A\otimes\ADir,A)\) is invertible
  (\cite{Meyer-Nest:BC}).  If \(b\in\KK^G_*(\C,A)\) and \(c\in\KK^G_*(A,\C)\),
  then we can write the Kasparov product \(b\otimes_A c\) as
  \[
  \C \xrightarrow{b} A
  \xleftarrow[\cong]{\ID_A\otimes\Dirac} A\otimes\ADir
  \xrightarrow{c\otimes\ID_\ADir} \ADir
  \xrightarrow{\Dirac} \C,
  \]
  where the arrows are morphisms in the category \(\KK^G\).  Therefore,
  \(b\otimes_A c\) factors through~\(\Dirac\) and hence is a boundary class by
  Theorem~\ref{the:Kreg}.
\end{proof}

\section{Dual-Dirac morphisms and the Carlsson--Pe\-der\-sen condition}
\label{sec:Carlsson-Pedersen}

Now we construct boundary classes in \(\RKK^G_*(\EG;\C,\C)\) from more
classical boundaries.  We suppose again that~\(\EG\) is \(G\)\nb-compact, so
that \(\EG\) is a coarse space.

A \emph{metrisable compactification} of \(\EG\) is a metrisable compact
space~\(Z\) with a homeomorphism between \(\EG\) and a dense open subset
of~\(Z\).  It is called \emph{coarse} if all scalar-valued functions on~\(Z\)
have vanishing variation; this implies the corresponding assertion for
operator-valued functions because \(C(Z,D)\cong C(Z)\otimes D\).
Equivalently, the embedding \(\EG\to Z\) factors through the Higson
compactification of~\(\EG\).  A compactification is called
\emph{\(G\)\nb-equivariant} if~\(Z\) is a \(G\)\nb-space and the embedding
\(\EG\to Z\) is \(G\)\nb-equivariant.  An equivariant compactification is
called \emph{strongly contractible} if it is \(H\)\nb-equivariantly
contractible for all compact subgroups~\(H\) of~\(G\).  The
\emph{Carlsson--Pe\-der\-sen condition} requires that there should be a
\(G\)\nb-compact model for~\(\EG\) that has a coarse, strongly contractible,
and \(G\)\nb-equivariant compactification.

Typical examples of such compactifications are the Gromov boundary for a
hyperbolic group (viewed as a compactification of the Rips complex), or the
visibility boundary of a CAT(0) space on which~\(G\) acts properly,
isometrically, and cocompactly.

\begin{theorem}
  \label{the:dual_Dirac_compactify}%
  Let~\(G\) be a locally compact group with a \(G\)\nb-compact model for
  \(\EG\) and let \(\EG\subseteq Z\) be a coarse, strongly contractible,
  \(G\)\nb-equivariant compactification.  Then~\(G\) has a dual-Dirac
  morphism.
\end{theorem}

\begin{proof}
  We use the \(C^*\)\nb-algebra \(\bbbarred_G(Z)\) as defined in
  \textsection\ref{sec:Higson_corona}.  Since~\(Z\) is coarse, there is an
  embedding \(\bbbarred_G(Z)\subseteq\vvr_G(\EG)\).  Let \(\partial Z\defeq
  Z\setminus \EG\) be the boundary of the compactification.  Identifying
  \[
  \bbbarred_G(\partial Z)\cong\bbbarred_G(Z)/C_0(\EG,\Comp_G),
  \]
  we get a morphism of extensions
  \[
  \xymatrix{
    0 \ar[r] & C_0(\EG,\Comp_G) \ar[r] \ar@{=}[d] &
    \vvr_G(\EG) \ar[r] &
    \bvr_G(\EG) \ar[r] & 0 \\
    0 \ar[r] & C_0(\EG,\Comp_G) \ar[r] &
    \bbbarred_G(Z) \ar[r] \ar[u]_{\subseteq} &
    \bbbarred_G(\partial Z) \ar[r] \ar[u]_{\subseteq} & 0.
  }
  \]
  Let~\(H\) be a compact subgroup.  Since~\(Z\) is compact, we have
  \(\bbbar_G(Z)=C(Z,\Comp)\).  Since~\(Z\) is \(H\)\nb-equivariantly
  contractible by hypothesis, \(\bbbar_G(Z)\) is \(H\)\nb-equivariantly
  homotopy equivalent to~\(\C\).  Hence \(\bbbarred_G(Z)\) has vanishing
  \(H\)\nb-equivariant \(\K\)\nb-theory.  This implies
  \(\Ktop_*\bigl(G,\bbbarred_G(Z)\bigr)= 0\) by~\cite{Chabert-Echterhoff-Oyono:Going}, so that the
  connecting map
  \begin{equation}
    \label{pineapple}
    \Ktop_{*+1}\bigl(G,\bbbarred_G(\partial Z)\bigr)
    \to \Ktop_*\bigl(G,C_0(\EG)\bigr)
    \cong \K_*(C_0(\EG)\cross G)
  \end{equation}
  is an isomorphism.   This in turn implies that the connecting map
  \[
  \Ktop_{*+1}\bigl(G,\bvr_G(\EG)\bigr)
  \to \Ktop_*\bigl(G,C_0(\EG)\bigr)
  \]
  is surjective.  Thus we can lift \(1\in \RKK^G_0(\EG;\C,\C) \cong
  \K_0(C_0(\EG)\cross G)\) to
  \[
  \alpha\in \Ktop_1\bigl(G,\bvr_G(\EG)\bigr)
  \cong \Ktop_1\bigl(G,\bvr_G(\abs{G})\bigr).
  \]
  Then \(\Psi_*(\alpha)\in\KK^G_0(\C,\ADir)\) is the desired dual-Dirac
  morphism.
\end{proof}

The group \(\Ktop_*\bigl(G,\bbbarred_G(\partial Z)\bigr)\) that appears in the
above argument is a \emph{reduced} topological \(G\)\nb-equivariant
\(\K\)\nb-theory for~\(\partial Z\) and hence differs from
\(\Ktop_*\bigl(G,C(\partial Z)\bigr)\).  The relationship between these two
groups is analysed in~\cite{Emerson-Meyer:Euler}.


\begin{bibdiv}
\begin{biblist}
\bib{Abels:Slices}{article}{
  author={Abels, Herbert},
  title={Parallelizability of proper actions, global \(K\)\nobreakdash -slices and maximal compact subgroups},
  journal={Math. Ann.},
  volume={212},
  date={1974/75},
  pages={1--19},
  review={\MRref {0375264}{51 \#11460}},
}

\bib{Chabert-Echterhoff-Oyono:Going}{article}{
  author={Chabert, J{\'e}r{\^o}me},
  author={Echterhoff, Siegfried},
  author={Oyono-Oyono, Herv{\'e}},
  title={Going-down functors, the K\"unneth formula, and the Baum--Connes conjecture},
  journal={Geom. Funct. Anal.},
  volume={14},
  date={2004},
  number={3},
  pages={491--528},
  issn={1016-443X},
  review={\MRref {2100669}{2005h:19005}},
}

\bib{Connes-Gromov-Moscovici:Lipschitz}{article}{
  author={Connes, Alain},
  author={Gromov, Mikhail},
  author={Moscovici, Henri},
  title={Group cohomology with Lipschitz control and higher signatures},
  journal={Geom. Funct. Anal.},
  volume={3},
  date={1993},
  number={1},
  pages={1--78},
  issn={1016-443X},
  review={\MRref {1204787}{93m:19011}},
}

\bib{Dranishnikov:Lipschitz}{article}{
  author={Dranishnikov, Alexander N.},
  title={Lipschitz cohomology, Novikov's conjecture, and expanders},
  language={Russian, with Russian summary},
  journal={Tr. Mat. Inst. Steklova},
  volume={247},
  date={2004},
  number={Geom. Topol. i Teor. Mnozh.},
  pages={59--73},
  issn={0371-9685},
  translation={ journal={Proc. Steklov Inst. Math.}, date={2004}, number={4 (247)}, pages={50--63}, issn={0081-5438}, },
  review={\MRref {2168163}{2006i:57054}},
}

\bib{Emerson-Meyer:Euler}{article}{
  author={Emerson, Heath},
  author={Meyer, Ralf},
  title={Euler characteristics and Gysin sequences for group actions on boundaries},
  journal={Math. Ann.},
  volume={334},
  date={2006},
  number={4},
  pages={853--904},
  issn={0025-5831},
  review={\MRref {2209260}{2007b:19006}},
}

\bib{Emerson-Meyer:Dualizing}{article}{
  author={Emerson, Heath},
  author={Meyer, Ralf},
  title={Dualizing the coarse assembly map},
  journal={J. Inst. Math. Jussieu},
  volume={5},
  date={2006},
  number={2},
  pages={161--186},
  issn={1474-7480},
  review={\MRref {2225040}{2007f:19007}},
}

\bib{Emerson-Meyer:Descent}{article}{
  author={Emerson, Heath},
  author={Meyer, Ralf},
  title={A descent principle for the Dirac--dual-Dirac method},
  journal={Topology},
  volume={46},
  date={2007},
  number={2},
  pages={185--209},
  issn={0040-9383},
  review={\MRref {2313071}{}},
}

\bib{Ferry-Ranicki-Rosenberg:Novikov}{article}{
  author={Ferry, Steven C.},
  author={Ranicki, Andrew},
  author={Rosenberg, Jonathan},
  title={A history and survey of the Novikov conjecture},
  conference={ title={}, address={Oberwolfach}, date={1993}, },
  book={ series={London Math. Soc. Lecture Note Ser.}, volume={226}, publisher={Cambridge Univ. Press}, place={Cambridge}, },
  date={1995},
  pages={7--66},
  review={\MRref {1388295}{97f:57036}},
}

\bib{Higson:Bivariant}{article}{
  author={Higson, Nigel},
  title={Bivariant $K$\nobreakdash -theory and the Novikov conjecture},
  journal={Geom. Funct. Anal.},
  volume={10},
  date={2000},
  number={3},
  pages={563--581},
  issn={1016-443X},
  review={\MRref {1779613}{2001k:19009}},
}

\bib{Kasparov:Novikov}{article}{
  author={Kasparov, Gennadi G.},
  title={Equivariant \(KK\)-theory and the Novikov conjecture},
  journal={Invent. Math.},
  volume={91},
  date={1988},
  number={1},
  pages={147--201},
  issn={0020-9910},
  review={\MRref {918241}{88j:58123}},
}

\bib{Meyer-Nest:BC}{article}{
  author={Meyer, Ralf},
  author={Nest, Ryszard},
  title={The Baum--Connes conjecture via localisation of categories},
  journal={Topology},
  volume={45},
  date={2006},
  number={2},
  pages={209--259},
  issn={0040-9383},
  review={\MRref {2193334}{2006k:19013}},
}


\end{biblist}
\end{bibdiv}
\end{document}